\documentclass{elsart}
\usepackage{natbib}
\usepackage{amsthm, amsmath, amsfonts}

\begin{document}

\begin{frontmatter}

\title{A fast and accurate algorithm for solving Bernstein--Vandermonde linear systems}

\author{Ana Marco \thanksref{EM1}}
\author{, Jos\'e-Javier Mart{\'\i}nez \thanksref{EM2}}

\address{Departamento de Matem\'aticas, Universidad
de Alcal\'a,}

\address{Campus Universitario, 28871-Alcal\'a de Henares (Madrid), Spain}

\thanks[EM1]{E-mail: ana.marco@uah.es}

\thanks[EM2]{Corresponding author. E-mail: jjavier.martinez@uah.es}

\begin{abstract}

A fast and accurate algorithm for solving a Bernstein-Vandermonde
linear system is presented. The algorithm is derived by using
results related to the bidiagonal decomposition of the inverse of a
totally positive matrix by means of Neville elimination. The use of
explicit expressions for the determinants involved in the process
serves to make the algorithm both fast and accurate.

\bigskip


\end{abstract}

\end{frontmatter}

\section{Introduction}

The {\it Bernstein basis} for the space of algebraic polynomials of
degree less than or equal to $n$ is a widely used basis in Computer
Aided Geometric Design due to the good properties that it possesses
(see, for instance, [5, 9, 10, 11, 17]). However, the explicit
conversion between the Bernstein and the power basis is
exponentially ill-conditioned as the polynomial degree increases
[10]. For this reason, it is very important that when designing
algorithms for performing numerical computations with polynomials
expressed in Bernstein form, all the intermediate operations are
developed using this form only [2]. A paper which presents various
basic operations for polynomials in Bernstein form is [11]. In [2]
an algorithm for computing the greatest common divisor of two
polynomials in Bernstein form avoiding explicit basis conversion is
given.

Our aim in this paper is to develop a fast and accurate algorithm
for solving a linear system whose coefficient matrix is a
Bernstein-Vandermonde matrix, that is to say, a Vandermonde-like
matrix for the Bernstein polynomials. Taking into account that this
matrix is the coefficient matrix of the linear system associated
with a Lagrange interpolation problem in the Bernstein basis, this
result will allow us to perform a basic polynomial procedure, the
Lagrange interpolation [7], avoiding transformations between
Bernstein and power basis. In this way, general algorithms for
polynomials in Bernstein form involving an interpolation stage could
be designed using this form only.

Our algorithm will be based on the bidiagonal factorization of the
inverse of the Bernstein-Vandermonde matrix. Factorizations in terms
of bidiagonal matrices are very useful when working with Vandermonde
[3, 14, 16], Cauchy [4], Cauchy-Vandermonde [20, 21] and generalized
Vandermonde matrices [8].

Let us observe here that, of course, Bernstein-Vandermonde linear
systems can be solved by using standard algorithms such as Gaussian
or Neville elimination. However they are not fast and
 the solution provided by them will generally be less accurate since
 Bernstein-Vandermonde matrices are ill conditioned (see Table 1, where $n$
 is the degree of the Bernstein polynomials involved in the definition of the matrix, whose order is $n+1$) and
these algorithms (which do not take into account the structure of
the matrix) can suffer from subtractive cancellation.

\begin{table}[h]
\begin{center}
{\scriptsize \caption{Condition number of a Bernstein-Vandermonde
matrix with $x_i=\frac{i}{n+2}$.} \label{tab1}
\begin{tabular}{|c|c|c|c|c|c|c|c|}
\hline
$n$ & $10$ & $15$ & $20$ & $25$ & $30$ & $35$ & $40$ \\
\hline $\kappa_{\infty}(A_n)$ & $2.1e+04$ & $2.6e+06$ & $3.5e+08$ &
$4.7e+10$  & $6.6e+12$ & $9.0e+14$ & $1.3e+17$ \\
\hline
\end{tabular}}
\end{center}
\end{table}

On the other hand, it must be observed that the high condition
numbers are due to the high norm of the inverse matrix, since as it
will be easily seen after the definitions of Section 3 the
$\infty$-norm of a Bernstein-Vandermonde matrix is always 1.
Therefore, taking into account that if $Ax=b$ and $A(x+\Delta
x)=b+\Delta b$ then $\Delta x=A^{-1} \Delta b$, the fact that $\Vert
A^{-1} \Vert_{\infty}$ is high implies that the effect of
perturbations in the vector $b$ is likely to be important. A lot
information about the related concepts of perturbation theory and
numerical stability of algorithms in the context of solving linear
systems can be found in Chapter 7 of [16].

The rest of the paper is organized as follows. Neville elimination
and total positivity are considered in Section 2. In Section 3 the
bidiagonal factorization of the inverse of a Bernstein-Vandermonde
matrix is presented. In Section 4 the algorithm for solving a linear
system whose coefficient matrix is Bernstein-Vandermonde, and the
computation of its complexity are given. Finally, Section 5 is
devoted to illustrate the accuracy of the algorithm by means of some
numerical experiments.

\section{Basic results on Neville elimination and total positivity}

In this section we will briefly recall some basic results on Neville
elimination and total positivity which we will apply in Section 3.
Our notation follows the notation used in [12] and [13]. Given $k$,
$n \in {\bf N}$  ($1 \leq k \leq n$), $Q_{k,n}$ will denote the set
of all increasing sequences of $k$ natural numbers less than or
equal to $n$.

Let $A$ be a real square matrix of order $n$. For $k \leq n$, $m
\leq n$, and for any $\alpha \in Q_{k,n}$ and $\beta \in Q_{m,n}$,
we will denote by $A[\alpha \vert \beta ]$ the submatrix $k\times m$
of $A$ containing the rows numbered by $\alpha $ and the columns
numbered by $\beta $.

The fundamental tool for obtaining the results presented in this
paper is the {\it Neville elimination} [12, 13], a procedure that
makes zeros in a matrix adding to a given row an appropriate
multiple of the previous one. For a nonsingular matrix
$A=(a_{i,j})_{1\leq i,j\leq n}$, it consists on $n-1$ steps
resulting in a sequence of matrices $A:=A_1\to A_2\to \ldots \to
A_n$, where $A_t= (a_{i,j}^{(t)})_{1\leq i,j\leq n}$ has zeros below
its main diagonal in the $t-1$ first columns. The matrix $A_{t+1}$
is obtained from $A_t$ ($t=1,\ldots ,n$) by using the following
formula:
$$
a_{i,j}^{(t+1)}:= \left\{ \begin{array}{ll} a_{i,j}^{(t)} &
\text{if} \quad i\leq t\\
a_{i,j}^{(t)}-(a_{i,t}^{(t)}/a_{i-1,t}^{t})a_{i-1,j}^{(t)}~ &
\text{if} \quad i\geq t+1 ~\text{and}~ j\geq t+1\\
0 & \text{otherwise}
\end{array}\right..  \eqno(2.1)
$$
In this process the element
$$
p_{i,j}:=a_{i,j}^{(j)} \qquad 1\leq j\leq n; ~~ j\leq i\leq n
$$
is called {\it pivot} ($i,j$) of the Neville elimination of $A$. The
process would break down if any of the pivots $p_{i,j}$ ($j\leq
i<n$) is zero. In that case we can move the corresponding rows to
the bottom and proceed with the new matrix, as described in [12].
The Neville elimination can be done without row exchanges if all the
pivots are nonzero. The pivots $p_{i,i}$ are called {\it diagonal
pivots}. If all the pivots $p_{i,j}$ are nonzero, then
$p_{i,1}=a_{i,1}\, \forall i$ and, by Lemma 2.6 of [12]
$$p_{i,j}={\det A[i-j+1,\ldots ,i\vert 1,\ldots ,j]\over \det
A[i-j+1,\ldots ,i-1\vert 1,\ldots ,j-1]}\qquad 1<j\leq i\leq n.
\eqno(2.2)$$

The element
$$
m_{i,j}=\frac{p_{i,j}}{p_{i-1,j}} \qquad 1\leq j\leq n; ~~ j< i\leq
n  \eqno(2.3)
$$
is called {\it multiplier} of the Neville elimination of $A$. The
matrix $U:=A_n$ is upper triangular and has the diagonal pivots in
its main diagonal.

The {\it complete Neville elimination} of a matrix $A$ consists on
performing the Neville elimination of $A$ for obtaining $U$ and then
continue with the Neville elimination of $U^T$. The pivot
(respectively, multiplier) $(i,j)$ of the complete Neville
elimination of $A$ is the pivot (respectively, multiplier) $(j,i)$
of the Neville elimination of $U^T$, if $j\ge i$. When no row
exchanges are needed in the Neville elimination of $A$ and $U^T$, we
say that the complete Neville elimination of $A$ can be done without
row and column exchanges, and in this case the multipliers of the
complete Neville elimination of $A$ are the multipliers of the
Neville elimination of $A$ if $i\ge j$ and the multipliers of the
Neville elimination of $A^T$ if $j\ge i$.

A matrix is called {\it totally positive} (respectively, {\it
strictly totally positive}) if all its minors are nonnegative
(respectively, positive). The Neville elimination characterizes the
strictly totally positive matrices as follows [12]:

{\bf Theorem 2.1.} A matrix is strictly totally positive if and only
if its complete Neville elimination can be performed without row and
column exchanges, the multipliers of the Neville elimination of $A$
and $A^T$ are positive, and the diagonal pivots of the Neville
elimination of $A$ are positive.

It is well known [5] that the Bernstein-Vandermonde matrix is a
strictly totally positive matrix when the interpolation points
satisfy $0<x_1<x_2<\ldots<x_{n+1}<1$, and this fact has inspired our
search for a fast algorithm, but this result will also be shown to
be a consequence of our Theorem 3.3.

\section{Bidiagonal factorization}

The {\it Bernstein basis} of the space $\Pi_n(x)$ of polynomials of
degree less than or equal to $n$ on the interval $[0,1]$ is:
$$
\mathcal{B}_n=\big\{ b_i^{(n)}(x) = {n \choose i} (1 - x)^{n-i} x^i,
\qquad i = 0, \ldots, n \big\}.
$$

We will call the following Vandermonde-like matrix for the Bernstein
basis $\mathcal{B}_n$,
$$
A=
\begin{pmatrix}
{n \choose 0}(1-x_1)^n & {n \choose 1}x_1(1-x_1)^{n-1}& \cdots & {n
\choose n}x_1^n\\
{n \choose 0}(1-x_2)^n & {n \choose
1}x_2(1-x_2)^{n-1}& \cdots & {n \choose n}x_2^n\\
\vdots & \vdots & \ddots & \vdots\\ {n \choose 0}(1-x_{n+1})^n & {n
\choose 1}x_{n+1}(1-x_{n+1})^{n-1}& \cdots & {n \choose n}x_{n+1}^n
\end{pmatrix},
$$
a {\it Bernstein-Vandermonde matrix}. Let us observe that these
Vandermonde-like matrices are not exactly the Vandermonde-like
matrices considered in [16], since all the polynomials in the
Bernstein basis $\mathcal{B}_n$ have the same degree $n$ and they do
not satisfy a three-term recurrence relation. From now on, we will
assume $0<x_1<x_2< \ldots <x_{n+1}<1$.

This matrix $A$ is the coefficient matrix of the linear system
associated with the following Lagrange interpolation problem in the
Bernstein basis $\mathcal{B}_n$: given the interpolation nodes
$\{x_i: ~ i=1, \ldots, n+1\}$ with $0<x_1<x_2< \ldots <x_{n+1}<1$
and the interpolation data $\{b_i: ~ i=1, \ldots, n+1\}$ find the
polynomial
$$
p(x)=\sum_{k=0}^n a_k{n \choose k}(1-x)^{n-k}x^k
$$
such that $p(x_i)=b_i$ for $i=1, \ldots, n+1$.  A good introduction
to the interpolation theory can be seen in [7].

\medskip
{\bf Proposition 3.1.}
$$
\det A = {n \choose 0}{n \choose 1} \cdots {n \choose n}
\prod_{1\leq i < j \leq n+1}(x_j-x_i).
$$

{\bf Proof.} It is easy to see that the matrix of change of basis
from the Bernstein basis $\mathcal{B}_n$ to the power basis
$\{1,x,x^2,\ldots,x^n\}$ is a lower triangular matrix of order $n+1$
whose diagonal elements are ${n \choose 0},{n \choose 1}, \ldots, {n
\choose n}$.

From this fact, it is obtained that
$$
\det A={n \choose 0}{n \choose 1} \cdots {n \choose n} \det V,
$$
where $V$ is the Vandermonde matrix
$$
V=\begin{pmatrix}
1 & x_1 & x_1^2 & \cdots & x_1^n\\
1 & x_2 & x_2^2 & \cdots & x_2^n\\
\vdots & \vdots & \ddots & \vdots\\
1 & x_{n+1} & x_{n+1}^2 & \cdots & x_{n+1}^n
\end{pmatrix}.
$$
Using the well-known formula for the determinant of  a Vandermonde
matrix
$$
\det V=\prod_{1\leq i < j \leq n+1}(x_j-x_i)
$$
the proof is concluded. $\Box$

\medskip
The next corollary follows directly from Proposition 3.1., and will
be useful to make the derivation of the algorithm easier.

\medskip
{\bf Corollary 3.2.}
$$
\det
\begin{pmatrix}
(1-x_1)^n & x_1(1-x_1)^{n-1}& \cdots & x_1^n\\
(1-x_2)^n & x_2(1-x_2)^{n-1}& \cdots & x_2^n\\
\vdots & \vdots & \ddots & \vdots\\ (1-x_{n+1})^n &
x_{n+1}(1-x_{n+1})^{n-1}& \cdots & x_{n+1}^n
\end{pmatrix}
=\prod_{1\leq i < j \leq n+1}(x_j-x_i)
$$

The following result will be the key to construct our algorithm.

\medskip
{\bf Theorem 3.3.} Let $A=(a_{i,j})_{1\le i,j\le n+1}$ be a
Bernstein-Vandermonde matrix whose nodes satisfy $0 < x_1 < x_2
<\ldots < x_n <x_{n+1} <1$. Then $A^{-1}$ admits a factorization in
the form
$$A^{-1}=G_1G_2\cdots G_{n}D^{-1}F_{n}F_{n-1}\cdots F_1, \eqno(3.1)$$
where $G_i$ are upper triangular bidiagonal matrices, $F_i$ are
lower triangular bidiagonal matrices ($i=1,\ldots,n$), and $D$ is a
diagonal matrix.

{\bf Proof.} The matrix $A$ is strictly totally positive (see [5])
and therefore, by Theorem 2.1, the complete Neville elimination of
$A$ can be performed without row and column exchanges providing the
following factorization of $A^{-1}$ (see [12] and [13]):
$$A^{-1}=G_1G_2\cdots G_{n}D^{-1}F_{n}F_{n-1}\cdots F_1, $$
where $F_i$ ($1\le i\le n$) are bidiagonal matrices of the form
$$F_i=\begin{pmatrix}
1 & & & & & & & \\
0 & 1 & & & & & & \\
& \ddots & \ddots & & & & & \\
& & 0 & 1 & & & & \\
& & & -m_{i+1,i} & 1 & & & \\
& & & & -m_{i+2,i} & 1 & & \\
& & & & & \ddots & \ddots & \\
& & & & & & -m_{n+1,i} & 1
\end{pmatrix},  \eqno(3.2)
$$
$G^T_i$ ($1\le i\le n$) are bidiagonal matrices of the form
$$G_i^T=\begin{pmatrix}
1 & & & & & & & \\
0 & 1 & & & & & & \\
& \ddots & \ddots & & & & & \\
& & 0 & 1 & & & & \\
& & & -\widetilde m_{i+1,i} & 1 & & & \\
& & & & -\widetilde m_{i+2,i} & 1 & & \\
& & & & & \ddots & \ddots & \\
& & & & & & -\widetilde m_{n+1,i} & 1
\end{pmatrix},  \eqno(3.3)
$$
and $D$ is the diagonal matrix whose $i$th ($1\le i\le n+1$)
diagonal entry is the diagonal pivot $p_{i,i}=a_{i,i}^{(i)}$ of the
Neville elimination of $A$:
$$
D=\text{diag}\{p_{1,1},p_{2,2}\ldots,p_{n+1,n+1}\}. \eqno(3.4)
$$

Taking into account that the minors of $A$ with $j$ initial
consecutive columns and $j$ consecutive rows starting with row $i$
are
$$
\begin{aligned}
& \det A [i,\ldots, i+j-1 \vert 1, \ldots, j] ={n \choose 0} {n
\choose 1} \cdots {n \choose j-1}\\ & (1-x_i)^{n-j+1}
(1-x_{i+1})^{n-j+1} \cdots (1-x_{i+j-1})^{n-j+1} \prod_{i \leq k <
l \leq i+j-1}(x_l-x_k),
\end{aligned}
$$
a result that follows from the properties of the determinants and
Corollary 3.2, and that $m_{i,j}$ are the multipliers of the Neville
elimination of $A$, we obtain that
$$
m_{i,j}=\frac{p_{i,j}}{p_{i-1,j}}=\frac{ (1-x_i)^{n-j+1} (1-x_{i-j})
\prod_{k=1}^{j-1}(x_i-x_{i-k}) }{ (1-x_{i-1})^{n-j+2}
\prod_{k=2}^j(x_{i-1}-x_{i-k}) },  \eqno(3.5)
$$
where $j=1,\ldots,n$ and $i=j+1, \dots, n+1$.

As for the minors of $A^T$ with $j$ initial consecutive columns and
$j$ consecutive rows starting with row $i$, they are:
$$
\begin{aligned}
& \det A^T[i,\ldots, i+j-1 \vert 1, \ldots, j] ={n \choose i-1} {n
\choose i} \cdots {n \choose i+j-2}x_1^{i-1}x_2^{i-1} \cdots
x_j^{i-1}\\
& (1-x_1)^{n-i-j+2}(1-x_2)^{n-i-j+2} \cdots
(1-x_j)^{n-i-j+2}\prod_{1 \leq k < l \leq j}(x_l-x_k).
\end{aligned}
$$
This expression also follows from the properties of the determinants
and Corollary 3.2. Since the entries $\widetilde m_{i,j}$ are the
multipliers of the Neville elimination of $A^T$, using the previous
expression for the minors of $A^T$ with initial consecutive columns
and consecutive rows, it is obtained that
$$
\widetilde m_{i,j}=\frac{(n-i+2)\cdot x_j}{(i-1)(1-x_j)} \qquad
\quad j=1,\dots,n;
 ~ i=j+1,\ldots,n+1. \eqno(3.6)
$$
Finally, the $i$th diagonal element of $D$
$$
p_{i,i}=\frac{ {n \choose i-1}(1-x_i)^{n-i+1} \prod_{k<i}(x_i-x_k)
}{ \prod_{k=1}^{i-1}(1-x_k) } \qquad i=1,\ldots,n+1, \eqno(3.7)
$$
is obtained by using the expression for the minors of $A$ with
initial consecutive columns and initial consecutive rows. $\Box$

Moreover, by using the same arguments of [20], it can be seen that
this factorization is unique among factorizations of this type, that
is to say, factorizations in which the matrices involved have the
properties shown by formulae (3.2), (3.3) and (3.4).

\medskip
Let us observe that the formulae obtained in the proof of Theorem
3.3 for the minors of $A$ with $j$ initial consecutive columns and
$j$ consecutive rows, and for the minors of $A^T$ with $j$ initial
consecutive columns and $j$ consecutive rows show that they are not
zero and so, the complete Neville elimination of $A$ can be
performed without row and column exchanges. Looking at equations
(3.5)-(3.7) is easily seen that $m_{i,j}$, $\widetilde m_{i,j}$ and
$p_{i,i}$ are positive. Therefore, taking into account Theorem 2.1,
this confirms that the matrix $A$ is strictly totally positive.

\section{The algorithm}

In this section we will present a fast algorithm for solving a
linear system whose coefficient matrix is a Bernstein-Vandermonde
matrix. In order to solve the linear system $Ax=b$, where $A$ is the
$(n+1)\times(n+1)$ Bernstein-Vandermonde matrix introduced in
Section 3, we use Theorem 3.3 for obtaining
$$x=A^{-1}b=G_1G_2\cdots G_{n}D^{-1}F_{n}F_{n-1}\cdots F_1b.$$

Since $F_i$ and $G_i$ ($i=1,\ldots,n+1$) are bidiagonal matrices and
$D^{-1}$ is a diagonal matrix, it is clear that the computational
complexity of computing the whole product from right to left is
$O(n^2)$. It remains to see that the construction of the matrices
$F_i$, $G_i$ and $D^{-1}$ can be carried out with a computational
complexity of $O(n^2)$.

Let us start with the entries $m_{i,j}$ given by equation (3.5). We
will use the following expressions
$$
\begin{aligned}
M_{i,1} & =\frac{(1-x_i)^n}{(1-x_{i-1})^{n+1}}\\
m_{i,1} & = (1-x_{i-1})\cdot M_{i,1}\\
M_{i,j} & =  \frac{(1-x_i)^{n-j+1} \prod_{k=1}^{j-1}(x_i-x_{i-k})}{(1-x_{i-1})^{n-j+2} \prod_{k=2}^j(x_{i-1}-x_{i-k})}\\
M_{i,j+1} & =
\frac{(1-x_{i-1})(x_i-x_{i-j})}{(1-x_i)(x_{i-1}-x_{i-j-1})}\cdot
M_{i,j}\\
m_{i,j+1} & =(1-x_{i-j-1})\cdot M_{i,j+1},
\end{aligned}
$$
where $i=2,\ldots,n+1$ and $j=1,\ldots,i-2$, in their construction:

{\tt for} $i=2:n+1$

\quad $M=\frac{(1-x_i)^n}{(1-x_{i-1})^{n+1}}$

\quad $m_{i,1}=(1-x_{i-1})\cdot M$

\quad {\tt for} $j=1:i-2$

\quad \quad
$M=\frac{(1-x_{i-1})(x_i-x_{i-j})}{(1-x_i)(x_{i-1}-x_{i-j-1})}\cdot
M$

\quad \quad $m_{i,j+1}=(1-x_{i-j-1})\cdot M$

\quad {\tt end}

{\tt end}

Now we compute the entries $\widetilde m_{i,j}$ given by equation
(3.6):

{\tt for} $j=1:n$

\quad $c_j=\frac{x_j}{1-x_j}$

\quad {\tt for} $i=j+1:n+1$

\quad \quad $\widetilde m_{i,j}=\frac{n-i+2}{i-1}\cdot c_j$

\quad {\tt end}

{\tt end}

As for the diagonal elements $p_{i,i}$ of $D$ given by equation
(3.7), they are constructed using the equalities
$$
\begin{aligned}
q_{1,1} & =1\\
p_{1,1} & =(1-x_1)^n\\
q_{i,i} & = \frac{{n \choose i-1}}{\prod_{k=1}^{i-1}(1-x_k)}\\
q_{i+1,i+1} & =\frac{n-i+1}{i\cdot(1-x_i)}\cdot q_{i,i} \\
p_{i+1,i+1} & =q_{i+1,i+1} \cdot (1-x_{i+1})^{n-i}
\prod_{k<i+1}(x_{i+1}-x_k),
\end{aligned}
$$
where $i=1, \ldots,n$, in the following way:

$q=1$

$p_{1,1}=(1-x_1)^n$

{\tt for} $i=1:n$

\quad $q=\frac{(n-i+1)}{i(1-x_{i})}\cdot q$

\quad $aux=1$

\quad {\tt for} $k=1:i$

\quad \quad $aux=(x_{i+1}-x_k)\cdot aux$

\quad {\tt end}

\quad $p_{i+1,i+1}=q\cdot (1-x_{i+1})^{n-i} \cdot aux$

{\tt end}

Looking at this algorithm is enough to conclude that the
computational complexity of the construction of the matrices $D$,
$F_i$ and $G_i$ ($i=1,\ldots,n+1$) is $O(n^2)$, and therefore, the
computational complexity of solving the whole linear system is also
$O(n^2)$.

A similar algorithm with computational complexity $O(n^2)$ can be
developed for computing the bidiagonal factorization of the inverse
of $A^T$, that is, $A^{-T}$. In consequence, an algorithm with
computational complexity $O(n^2)$ is obtained for solving the dual
linear system $A^Tx=b$, where $A$ is an $(n+1)\times(n+1)$
Bernstein-Vandermonde matrix, by using:
$$
x=A^{-T}b=F_1^TF_2^T\cdots F_n^TD^{-1}G_n^TG_{n-1}^T\cdots G_1^T.
$$

\section{Numerical experiments and final remarks}

Finally we present some numerical experiments which illustrate the
good properties of our algorithm. We compute the exact solution
$x_e$ of each one of the Bernstein-Vandermonde linear systems $Ax=b$
by using the command {\tt linsolve} of {\it Maple 10} and use it for
comparing the accuracy of the results obtained in M{\scriptsize
ATLAB} by means of:
\begin{enumerate}

\item The algorithm presented in Section 4 for computing the
bidiagonal decomposition of $A^{-1}$. We will call it {\tt MM}.

\item The algorithm {\tt TNBD} of Plamen Koev [18] that
computes the bidiagonal decomposition of $A^{-1}$ without taking
into account the structure of $A$.

\item The command $A\backslash b$ of M{\scriptsize ATLAB}.

\end{enumerate}

The fast product (from right to left) of the bidiagonal matrices and
the vector $b$ is also implemented in M{\scriptsize ATLAB} and is
the second stage in the solution of the linear system in (1) and
(2).

We compute the relative error of a solution $x$ of the linear system
$Ax=b$ by means of the formula:
$$
err=\frac{\Vert x - x_e \Vert_2}{\Vert x_e \Vert_2}.
$$

{\it Remark.} The algorithm {\tt TNBD} computes the matrix denoted
as $\mathcal{BD}(A)$ in [19], which represents the {\it bidiagonal
decomposition} of $A$. But it is a remarkable fact that the same
matrix $\mathcal{BD}(A)$ also serves to represent the bidiagonal
decomposition of $A^{-1}$. The algorithm computes $\mathcal{BD}(A)$
by performing {\it Neville elimination} on $A$, which involves true
substractions, and therefore does not guarantee high relative
accuracy.

A detailed error analysis of Neville elimination, which shows the
advantages of this type of elimination for the class of totally
positive matrices, has been carried out in [1], and related work for
the case of Vandermonde linear systems can be seen in Chapter $22$
of [16].

{\bf Example 5.1} Let $\mathcal{B}_{10}$ be the Bernstein basis of
the space of polynomials with degree less than or equal to $10$ on
$[0,1]$ and $A$ be the Bernstein-Vandermonde matrix of order $11$
generated by the nodes $x_i=\frac{i}{12}$ for $i=1,\ldots,11$. The
condition number of $A$ is $\kappa_2(A)=1.8e+04$. Let us consider
$$
\begin{aligned}
b_1^T & =(1,0,2,-1,3,1,-2,0,0,3,5)^T\\
b_2^T & =(1,-2,1,-1,3,-1,2,-1,4,-1,1)^T
\end{aligned}
$$
two vectors of data.

The relative errors obtained when using the approaches (1), (2) and
(3) for solving the systems $Ax=b_i$ ($i=1,2$) are reported in Table
2.

\begin{table}[h]
\begin{center}
\begin{tabular}{|c|c|c|c|}
\hline $b_i$  & $MM$ & $TNBD$ & $A\backslash b_i$ \\
\hline $b_1$ & 1.3e-15 & 7.8e-14 &  5.4e-14\\
\hline $b_2$ & 8.6e-16 & 8.2e-14 &  1.0e-14\\
\hline
\end{tabular}
\end{center}\caption{Example 5.1: $A$ is a Bernstein-Vandermonde matrix of order
11.}
\end{table}

The following example will show how the accuracy of our approach is
maintained when the order, and therefore the condition number (see
Table 1), of the Bernstein-Vandermonde matrix increases, while the
accuracy of the other two approaches which do not exploit the
structure of the matrix $A$ goes down.

{\bf Example 5.2} Let $\mathcal{B}_{15}$ be the Bernstein basis of
the space of polynomials with degree less than or equal to $15$ on
$[0,1]$ and $A$ be the Bernstein-Vandermonde matrix of order $16$
generated by the nodes $x_i=\frac{i}{17}$ for $i=1,\ldots,16$. The
condition number of $A$ is $\kappa_2(A)=2.3e+06$. Let
$$
\begin{aligned}
b_1^T & =(2, 1, 2, 3, -1, 0, 1, -2, 4, 1, 1, -3, 0, -1, -1, 2)^T \\
b_2^T & =(1, -2, 1, -1, 3, -1, 2, -1, 4, -1, 2, -1, 1, -3, 1, -4)^T
\end{aligned}
$$
two vectors containing the data. The relative errors of the
solutions of the linear systems $Ax=b_i$ ($i=1,2$) obtained by means
of the approaches (1), (2) and (3) are reported in Table 3.

\begin{table}[h]
\begin{center}
\begin{tabular}{|c|c|c|c|}
\hline $b_i$  & $MM$ & $TNBD$ & $A\backslash b_i$ \\
\hline $b_1$ & 1.0e-15 & 5.9e-11 & 6.5e-12 \\
\hline $b_2$ & 4.9e-16 & 5.9e-11 & 6.4e-12 \\
\hline
\end{tabular}
\end{center}\caption{Example 5.2: $A$ is a Bernstein-Vandermonde matrix of order 16.}
\end{table}

Let us observe that in the first stage of the algorithm, which
corresponds to the computation of the bidiagonal decomposition of
$A^{-1}$, the {\it high relative accuracy} of the algorithm is
obtained because no subtractive cancellation occurs: we are
multiplying, dividing, or adding quantities with the same sign, or
forming $1-x_i$ and $x_i\pm x_j$ where $x_i$ and $x_j$ are initial
data.

As for the second stage of the algorithm corresponding to the
evaluation of the product
$$
G_1G_2\cdots G_{n}D^{-1}F_{n}F_{n-1}\cdots F_1b,
$$
it must be noted that when the vector $b$ has alternating sign
pattern ( $sign(b_i)=(\pm1)^i$), this stage is substraction-free,
and so the product can be computed with high relative accuracy. This
is a consequence to the checkerboard sign pattern of $A^{-1}$, which
derives from the fact that $A$ is a totally positive matrix.

This important property was already observed in an analogous
situation in the paper [15], devoted to the error analysis of the
Bj\"orck-Pereyra algorithm for Vandermonde systems. There the use of
the increasing ordering for the interpolation points is recommended,
an ordering which in the case $0<x_1<x_2<\dots<x_n$ makes the
Vandermonde matrix a totally positive one.

Our next two examples will serve to illustrate the concept of {\it
effective well-conditioning} introduced by Chan and Foulser in [6].
This concept has also been studied in the context of totally
positive Cauchy systems in [4], where it is seen that the {\it BKO}
algorithm exploits the effective well-conditioning to produce higher
accuracy for special right-hand sides.

Let $A \in R^{n \times n}$ with singular value decomposition
$A=U\Sigma V^T$, where $\Sigma=\text{diag}(\sigma_i)$, $\sigma_1
\geq \sigma_2 \geq \cdots \geq \sigma_n$ and $U=[u_1 u_2 \cdots
u_n]$. Let $P_k=U_kU_k^T$, with $U_k=[u_{n+1-k} \cdots u_n]$, the
projection operator onto the linear span of the smallest $k$ left
singular vectors of $A$. The {\it Chan-Foulser number} [4, 6] for
the linear system $Ax=f$ is defined by
$$
\gamma(A,f)=\min_k \frac{\sigma_{n-k+1}}{\sigma_n}\frac{\Vert f
\Vert_2}{\Vert P_kf \Vert_2}.
$$

Let us observe that the computation of the Chan-Foulser number
simplifies when $f=u_i$ ($i=1,\ldots,n$). In this case, using the
fact that the left singular vectors $u_1, \ldots, u_n$ are the
columns of the orthogonal matrix $U$, we obtain
$$
\gamma(A,u_i)=\frac{\sigma_{i}}{\sigma_n}.
$$

{\bf Example 5.3} Let $A$ be the same matrix as the one considered
in Example 5.2, whose condition number is $\kappa_2(A)=2.3e+06$. We
have solved the $16$ linear systems $Ax=u_i$, where $u_1, \ldots,
u_{16}$ are the left singular vectors of $A$. The results obtained
when using the approaches (1), (2) and (3), and the corresponding
Chan-Foulser numbers are reported in Table 4.

\begin{table}[h]
\begin{center}
\begin{tabular}{|c|c|c|c|c|}
\hline $u_i$ & $\gamma(A,u_i)$ & $MM$ & $TNBD$ & $A\backslash u_i$
\\
\hline $u_1$ & 2.3e+06 & 1.1e-10 & 4.7e-12 &  2.8e-11\\
\hline $u_2$ & 2.1e+06 & 5.0e-11 & 1.0e-11 &  1.7e-11\\
\hline $u_3$ & 1.7e+06 & 2.5e-11 & 1.3e-11 &  1.6e-10\\
\hline $u_4$ & 1.3e+06 & 4.9e-11 & 3.8e-11 &  3.4e-11\\
\hline $u_5$ & 9.3e+05 & 4.3e-11 & 1.1e-11 &  1.4e-11\\
\hline $u_6$ & 6.0e+05 & 3.1e-11 & 1.7e-11 &  6.2e-12\\
\hline $u_7$ & 3.5e+05 & 4.0e-11 & 7.5e-12 &  2.1e-11\\
\hline $u_8$ & 1.8e+05 & 1.8e-12 & 1.8e-11 &  1.4e-12\\
\hline $u_9$ & 8.5e+04 & 1.2e-11 & 9.8e-12 &  6.0e-12\\
\hline $u_{10}$ & 3.4e+04 & 1.7e-12 & 1.5e-11 &  9.1e-12\\
\hline $u_{11}$ & 1.2e+04 & 4.9e-13 & 2.0e-11 & 3.1e-12\\
\hline $u_{12}$ & 3.4e+03 & 6.5e-13 & 1.4e-11 &  1.4e-11\\
\hline $u_{13}$ & 7.9e+02 & 1.4e-13 & 2.4e-11 &  2.4e-11\\
\hline $u_{14}$ & 1.4e+02 & 8.1e-14 & 6.3e-12 &  2.4e-11\\
\hline $u_{15}$ & 1.6e+01 & 7.1e-15 & 2.1e-11 &  1.3e-11\\
\hline $u_{16}$ & 1 & 5.1e-16 & 5.9e-11 &  6.3e-12\\
\hline
\end{tabular}
\end{center}\caption{Example 5.3: left singular vectors $u_i$ for the right hand
side.}
\end{table}

{\bf Example 5.4} Let $\mathcal{B}_{15}$ be the Bernstein basis of
the space of polynomials with degree less than or equal to $15$ on
$[0,1]$ and $A$ be the Bernstein-Vandermonde matrix of order $16$
generated by the interpolation nodes

{\scriptsize $$
 \frac{1}{18} < \frac{1}{16} < \frac{1}{14} <
\frac{1}{12} < \frac{1}{10} < \frac{1}{8} < \frac{1}{6} <
\frac{1}{4} < \frac{11}{20} < \frac{19}{34} < \frac{17}{30} <
\frac{15}{26} < \frac{11}{18} < \frac{9}{14} < \frac{7}{10} <
\frac{5}{6}.
$$}
Its condition number is $\kappa_2(A)=3.5e+09$. We have solved the
$16$ linear systems $Ax=u_i$, where $u_1, \ldots, u_{16}$ are the
left singular vectors of $A$.

The results obtained when using the approaches (1), (2) and (3), and
the corresponding Chan-Foulser numbers are reported in Table 5.

\begin{table}
\begin{center}
\begin{tabular}{|c|c|c|c|c|}
\hline $u_i$ & $\gamma(A,u_i)$ & $MM$ & $TNBD$ & $A\backslash u_i$ \\
\hline $u_1$ & 3.5e+09 & 3.5e-07 & 3.2e-08 &  4.6e-08\\
\hline $u_2$ & 2.6e+09 & 1.1e-07 & 2.4e-08 &  7.3e-08\\
\hline $u_3$ & 1.3e+09 & 2.9e-08 & 1.3e-08 &  3.7e-08\\
\hline $u_4$ & 1.2e+09 & 2.2e-08 & 1.2e-08 &  5.0e-09\\
\hline $u_5$ & 5.3e+08 & 2.7e-08 & 2.3e-08 &  5.5e-09\\
\hline $u_6$ & 4.0e+08 & 2.5e-09 & 1.1e-08 &  4.6e-08\\
\hline $u_7$ & 1.1e+08 & 3.6e-09 & 2.4e-08 &  8.5e-09\\
\hline $u_8$ & 5.8e+07 & 2.6e-09 & 3.2e-09 &  8.0e-09\\
\hline $u_9$ & 1.1e+07 & 2.4e-10 & 1.1e-08 &  2.5e-08\\
\hline $u_{10}$ & 3.7e+06 & 3.9e-10 & 5.3e-09 &  8.1e-09\\
\hline $u_{11}$ & 4.8e+05 & 4.8e-12 & 8.0e-09 &  1.3e-08\\
\hline $u_{12}$ & 1.2e+05 & 1.5e-11 & 1.1e-08 &  2.2e-08\\
\hline $u_{13}$ & 6.2e+03 & 2.4e-12 & 1.5e-08 &  3.3e-08\\
\hline $u_{14}$ & 9.3e+02 & 7.3e-13 & 2.8e-10 &  1.9e-08\\
\hline $u_{15}$ & 1.4e+01 & 4.3e-14 & 3.0e-09 &  3.6e-08\\
\hline $u_{16}$ & 1 & 7.6e-15 & 1.1e-09 &  1.4e-08\\
\hline
\end{tabular}
\end{center}\caption{Example 5.4: left singular vectors $u_i$ for the right hand
side.}
\end{table}

The results appearing in Table 4 and Table 5 show that, as it could
be expected, the solution vectors obtained by conventional methods
based on Neville elimination or on Gaussian elimination suffer from
a loss of digits of precision which is roughly proportional to the
decimal logarithm of the condition number of the
Bernstein-Vandermonde matrix.

On the contrary, the results in Table 4 and Table 5 indicate that
the accuracy obtained by our algorithm, which exploits the structure
of the matrix, is governed by the Chan-Foulser number. In
particular, the very high relative precision observed in the 16th
system of both examples reflects two facts: the Chan-Foulser number
is equal to 1 and the vector $u_{16}$ has alternating sign pattern (
$sign(b_i)=(\pm1)^i$).

{\bf Acknowledgements}

This research has been partially supported by Spanish Research Grant
BFM 2003-03510 from the Spanish Ministerio de Ciencia y
Tecnolog{\'\i}a. We are grateful to the referees for their
suggestion of extending the numerical experiments.

\end{document}